\RequirePackage{fix-cm}
\documentclass[smallextended]{svjour3}       
\smartqed  

\usepackage{graphicx}
\usepackage{hyperref}
\usepackage{amsmath}
\usepackage[section]{placeins}
\graphicspath{ {./images/} }
\begin{document}

\title{Quadrature by Two Expansions for Evaluating Helmholtz Layer Potentials 
}

\titlerunning{QB2X for Helmholtz Layer Potentials}        

\author{Jared Weed \and Lingyun Ding \and Jingfang Huang \and Min Hyung Cho }


\institute{J. Weed \at
              Department of Mathematical Sciences \\
              University of Massachusetts Lowell\\
              \email{Jared\_Weed@student.uml.edu}           
           \and
           L. Ding \at
              Department of Mathematics \\
	     University of North Carolina at Chapel Hill\\
              \email{dingly@live.unc.edu}
	\and
	      J. Huang \at
              Department of Mathematics \\
	     University of North Carolina at Chapel Hill\\
              \email{huang@email.unc.edu}
              \and
              M.H. Cho \at
	      Corresponding Author \\
              Department of Mathematical Sciences \\
              University of Massachusetts Lowell\\
	      Tel.: 978-934-2410\\
              Fax: 978-934-3053 \\
              \email{minhyung\_cho@uml.edu}  
}

\date{Received: date / Accepted: date}

\maketitle

\begin{abstract}
In this paper, a Quadrature by Two Expansions (QB2X) numerical integration technique is developed 
for the single and double layer potentials of the Helmholtz equation in two dimensions. 
The QB2X method uses both local complex Taylor expansions and plane wave type expansions to achieve a
resulting representation which is numerically accurate for all target points (interior, exterior, or exactly
on the boundary) inside a leaf box in the fast multipole method (FMM) hierarchical tree structure. 
Compared to the original Quadrature by Expansion (QBX) method, the QB2X method explicitly
includes the nonlinearity from the boundary geometry in the plane wave expansions, thereby providing for
higher order representations of both the boundary geometry and density functions 
in the integrand, with its convergence following standard FMM error analysis. Numerical results 
are presented to demonstrate the performance of the QB2X method
for Helmholtz layer potentials and its comparison with the original QBX method for both flat and 
curved boundaries with various densities. The QB2X method overcomes the challenges of the 
original QBX method, and is better suited for efficient solutions of the Helmholtz equation 
with complex geometries.

\end{abstract}

\section{Introduction}
\label{intro}
Applying classical potential theory of partial differential equation (PDE) analysis, 
the solution for the homogeneous Helmholtz equation
\begin{equation}\begin{aligned}
\Delta u+k^2u=0
\end{aligned}\end{equation}
with a given boundary condition on the boundary $\Gamma$ can be represented as the linear combinations 
of single and double layer potentials 
\begin{equation}\begin{aligned}
S\psi(w) &= \int_{\Gamma} G(w,z) \psi(z) \mathrm{d}z \mbox{ and }
D\psi(w)= \int_{\Gamma} \frac{\partial G(w,z)}{\partial \mathbf{n}_z} \psi(z) \mathrm{d}z,
\end{aligned}\end{equation}
where $k$ is the wave number, $G(w, z)= \frac{\mathrm{i}}{4} H_0^{(1)} (k |w-z|)$ is the free-space Green's 
function at the target point $w$ due to the source point $z$ on the boundary, $H_0^{(1)}$ denotes 
the 0th order Hankel function of the first kind, $\mathbf{n}_z$ is the outward normal vector 
at $z$, and $\psi$ is the unknown density function (may be different for single and double layer potentials) 
on the boundary. To satisfy the given boundary condition, a common practice with integral equation methods 
is to find an appropriate combination of the layer potentials and derive a well-conditioned 
Fredholm second-kind boundary integral equation for the unknown density functions. Then, proper quadrature 
rules are applied to discretize the continuous integral equation formulation and yield a system of 
linear algebraic equations which can be solved by either fast direct or iterative solvers. 
One of the main numerical challenges with this process is the accurate and efficient 
evaluation of layer potentials with singular or nearly-singular kernels, especially when the 
target point $w$ is close to, or on the boundary of $\Gamma$. Note that the target point is located
exactly on the boundary in the boundary integral equation formulation, and in applications such 
as near-field optics \cite{lewis2003near}, surface plasmon resonance \cite{khattak2019linking}, 
and meta-materials \cite{tsantili2018computational}, the solution near the boundary plays an important 
role in determining the physical properties of the system. Numerical integration quadrature 
rules for singular and nearly-singular integrals have been extensively studied in literature. 
One technique is to modify the classical Newton-Cotes or Gaussian quadrature rules, either locally or 
globally, to compensate for the singularity in the kernel. The resulting schemes include 
an Alpert quadrature \cite{alpert1999hybrid}, a Rokhlin-Kapur quadrature \cite{kapur1997high}, 
a zeta correction quadrature \cite{wu2021zeta}, and a generalized Gaussian quadrature 
\cite{bremer2010nonlinear,yarvin1998generalized}. Other techniques involve analytically removing 
singularities using a change of variables \cite{bruno2001fast,duffy1982quadrature} or a  
regularization-correction method \cite{beale2001method}. Some of these quadrature techniques 
only work when the target point is located exactly on the boundary, and the performance of 
most existing quadrature rules depends on the complexity of the boundary geometry and distance 
between the target point and boundary $\Gamma$, which determines the kernel's singularity features.

\begin{figure}[b] 
   \centering
   \includegraphics[width=3.0in]{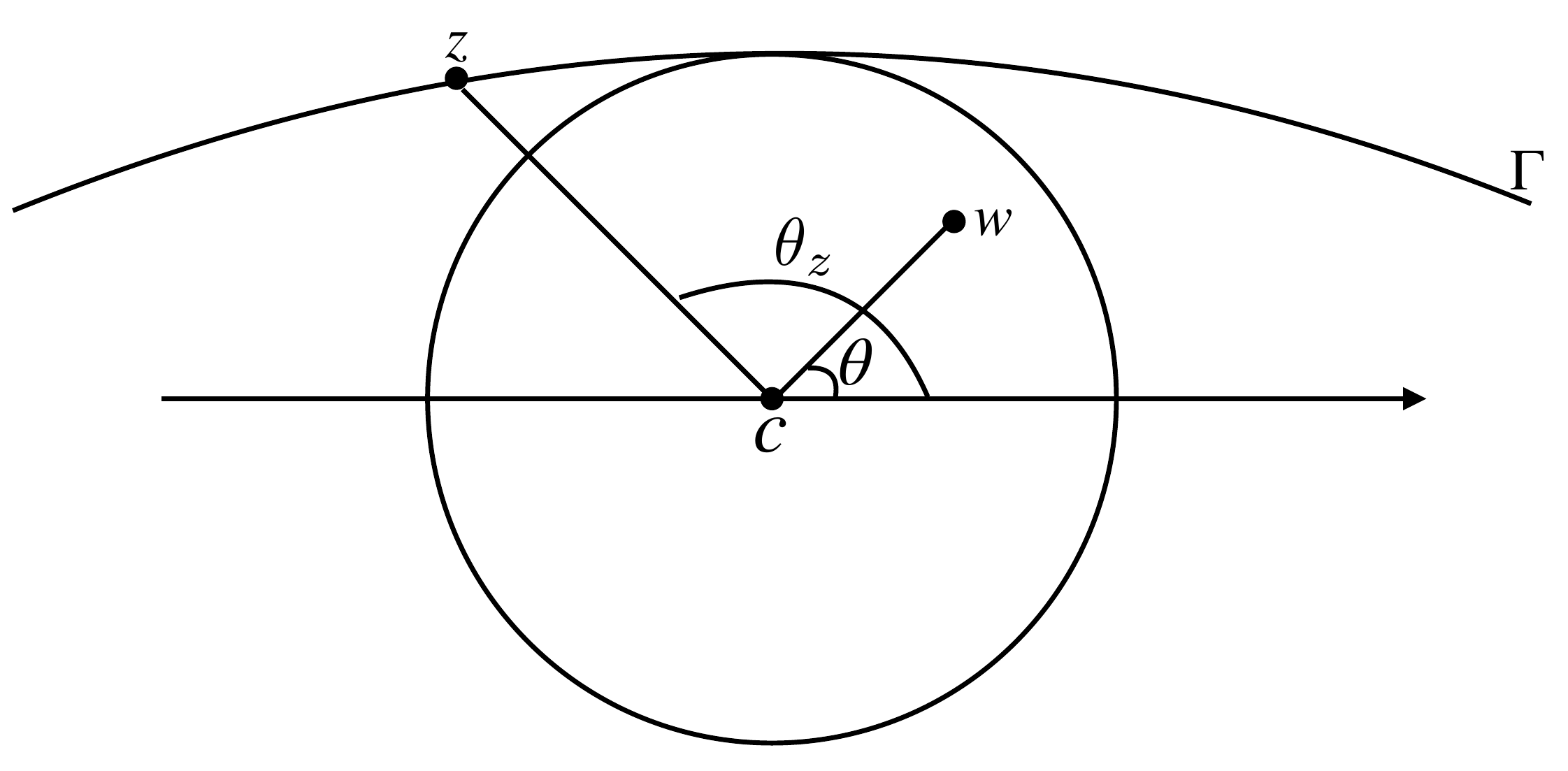} 
   \caption{Notation for the QBX method: $\Gamma$ is the boundary, $c$ is the center of expansion, 
	$w$ is the target point, $z$ is a point on $\Gamma$, $\theta_z$ is the angle of $z-c$ 
	measured from the $x$-axis, and $\theta$ is the angle of $w-c$ measured from the $x$-axis.}
   \label{fig:QBXnotation}
\end{figure}

By utilizing the smoothness properties of the layer potentials respectively in the interior 
and exterior of the domain, along with special partial wave or local expansions which automatically 
satisfy the underlying partial differential equations, recent research efforts have developed
non-traditional quadrature rules which are valid in a region (including 
the boundary) of the computational domain. One technique is the 
Quadrature by Expansion (QBX) method \cite{klockner2013quadrature,rachh2017fast}, 
which maps a layer potential to a partial wave (local) expansion centered at a point $c$ 
away from the boundary
\begin{equation}\begin{aligned}
\phi(w) = \sum_{n=-\infty}^{\infty} \alpha_n J_n(k||w-c||) e^{-\mathrm{i}n\theta}\label{eq:partial},
\end{aligned}\end{equation}
where
\begin{align}
&\alpha_n = \frac{\mathrm{i}}{4} \int_{\Gamma} H_n^{(1)}(k |z-c|)e^{\mathrm{i} n \theta_z}\psi(z) \mathrm{d}z \mbox{ for the Single layer potentials}, \nonumber\\
&\alpha_n = \frac{\mathrm{i}}{4} \int_{\Gamma} \frac{\partial}{\partial \mathbf{n}_z}H_n^{(1)}(k |z-c|)e^{\mathrm{i} n \theta_z}\psi(z) \mathrm{d}z \mbox{ for the Double layer potentials},\nonumber
\end{align}
$\theta$ and $\theta_z$ are the polar angles of $w-c$ and $z-c$, respectively 
(See Fig. \ref{fig:QBXnotation}), and $J_n$ and $H_n^{(1)}$ are the nth order Bessel and 
Hankel functions of the first kind, respectively. 
The value of the layer potential at a point reasonably close to $c$ is then approximated by 
the partial wave expansion. The partial wave expansion in Eq. (\ref{eq:partial}) is 
a direct consequence of applying Graf's Addition Theorem \cite{abramowitz1964handbook} to the 
Hankel function, and is the foundation of the fast multipole method (FMM) for the Helmholtz 
equation \cite{cho2010wideband,rokhlin1990rapid}. 

Existing numerical experiments show that when the boundary geometry is complex, 
the partial wave expansions in the QBX method may converge slowly due to the complicated nonlinear effects 
from the boundary geometry. The main contributions of this paper are new representations of the 
Helmholtz layer potentials that use both complex local Taylor expansions and plane wave type expansions
\begin{equation}\begin{aligned}
\phi(w) = \sum_{n=0}^N c_n (w-c)^n-2\pi \mathrm{i} \left( \sum_{p=0}^P  \sum_{j=1}^{J_u} \frac{b_p e^{\frac{\mathrm{i}p\pi}{2}r^u_j}}{1+\mathrm{i}s'(r^u_j)}-\sum_{p=-P}^{-1}  \sum_{j=1}^{J_d} \frac{b_p e^{\frac{\mathrm{i}p\pi}{2}r^d_j}}{1+\mathrm{i}s'(r^d_j)}\right),
\end{aligned}
\label{eq:QB2XIntro}
\end{equation}
where the boundary $\Gamma$ is defined by the parametric curve $z(t)=t+is(t)$, $s(t)$ is a real-coefficient 
polynomial for $t\in [-1 ,1]$, $w= x+\mathrm{i}y$ is a target point, $c$ is the center of the complex Taylor expansion, 
$c_n$ and $b_p$ are the complex coefficients of the Taylor expansion and plane wave type expansion, respectively, 
and $\{r_j^u\}_{j=1}^{J_u}$ and $\{r_j^d\}_{j=1}^{J_d}$ are the roots of the polynomial equation $w-(z+\mathrm{i}s(z))=0$ 
in the upper and lower complex half-plane, respectively. We refer to this new technique as the Quadrature 
by Two Expansions (QB2X). In linear algebra and harmonics analysis, the redundant basis sets (polynomial basis and 
exponential functions are linearly dependent) form a frame \cite{adcock2019frames,casazza2012finite,duffin1952class}.   
Compared with the QBX approach, our analysis and numerical results show that the convergence of the two expansions 
in QB2X follows standard FMM error analysis and the number of QB2X expansion terms can be explicitly 
determined by the degree of the polynomial $s(t)$ and number of terms in the Fourier extension approximation 
of the density function. Numerical results in Sec.~\ref{sec:3} show that for high order discretization 
of the boundary and oscillatory density functions, the QB2X method for Helmholtz layer potentials (QB2X-Helmholtz) 
converges rapidly and stably in the entire leaf box in the FMM tree structure.

In Eq.~(\ref{eq:QB2XIntro}), instead of a partial wave expansion for the Helmholtz equation, 
we use the polynomial expansion $\sum_{n=0}^N c_n (w-c)^n$ which only satisfies the Laplace 
equation for $w$. There are several reasons for this particular choice of basis function.
First, this allows the direct application of existing QB2X results for the Laplace 
equation (QB2X-Laplace) \cite{ding2021quadrature}. One can extract the logarithmic singularity from the 
Helmholtz Green's function~\cite{colton1998inverse} and apply the same contour integration 
and Residue Theorem techniques for the singular part. The remaining smooth (and analytical) 
part can be combined with the Fourier extension representation of the density function. 
Second, the numerical evaluation of a polynomial function is more efficient than that of 
the special Bessel or Hankel function. Note that the FMM handles the far field layer 
potentials using the multipole and local expansions based on the low-rank properties derived 
from the separation of variables of the kernel function. The near field singular or 
near singular layer potential evaluation is therefore only required for the low- to 
mild-frequency regimes when following the rule of ``10 points per wavelength" in the numerical
discretization. In this case, we expect there is no significant difference between 
the numbers of terms in the polynomial and partial wave expansions.  

We organize this paper as follows: In Sec. \ref{sec:2}, we present the derivation 
of the QB2X representation for the Helmholtz single layer potential. We omit the details
for the double layer potential as the analysis is nearly identical to that of the single layer case. 
In Sec. \ref{sec:3}, numerical results are presented to demonstrate the effectiveness of 
the QB2X method and its comparison with the QBX method. We conclude the paper in Sec. \ref{sec:4} 
with a summary and discussions of future work.  

\section{Quadrature by Two Expansions for Helmholtz Layer Potentials}
\label{sec:2}
In this section, we first show the parametrization of the single and double layer potentials, and analytically 
extract the logarithmic singularity from the Helmholtz kernel to allow for the easy adoption of existing 
QB2X-Laplace analysis results. Next, the QB2X-Helmholtz representations are derived for both straight line and 
curved boundaries. Finally, we present a stable numerical scheme for the case when the roots of the denominator are close 
to each other. 

\subsection{Parametrization of Helmholtz layer potentials}
Let the boundary $\Gamma$ be parametrized by $z(t)=t+\mathrm{i}s(t)$, for $-1 \leq t \leq 1$, where $s(t)$ is a real-coefficient polynomial
satisfying $s(0)=0$, $s'(0)=0$ after the proper rotations and translations. The single and double layer potentials at 
the target point $w=x+\mathrm{i}y$ are then parametrized as
\begin{equation}\begin{aligned}
S\psi(w) &= \int_{-1}^1 M(w,t) \psi(t) dt \mbox{ and } D\psi(w) =\int_{-1}^1 L(w,t) \psi(t) dt, \label{eq:para1}
\end{aligned}\end{equation}
where
\begin{equation}\begin{aligned}
M(w,t)&= \frac{i}{4}H_0^{(1)}(kr(w,t))|z'(t)|, \\
L(w,t)&= -\frac{ik}{4}\frac{(x-t)s'(t)-(y-s(t))}{r(w,t)}H^{(1)}_1(kr(w,t)), \\
r(w, t) &= \sqrt{(x-t)^2+(y-s(t))^2},\\
|z'(t)| &=\sqrt{1+s'(t)^2}.
\end{aligned}\end{equation}
Since evaluation of layer potentials are only required in the near-field direct interactions between the
leaf node and its neighbors (including itself) in the FMM algorithm, we apply the series expansions of the 
second kind Bessel function for small- and mild-frequency $kr(w,t)$ values (See \cite{colton1998inverse,kress2013linear}), 
and explicitly extract the logarithmic singularity of $M(w,t)$ and $L(w,t)$ as 
\begin{equation}\begin{aligned}
M(w,t) &=M_1(w,t)\ln{((x-t)^2+(y-s(t))^2)}+M_2(w,t), \\
L(w,t) &=L_1(w,t)\ln{((x-t)^2+(y-s(t))^2)}+L_2(w,t), 
\end{aligned}\end{equation}
where
\begin{equation}\begin{aligned}
M_1(w, t) &= -\frac{1}{4\pi}J_0(kr(w,t))|z'(t)|,\\
M_2(w, t) &= M(w,t)-M_1(w,t)\ln{((x-t)^2+(y-s(t))^2)},\\
L_1(w, t) &= \frac{k}{4\pi}((x-t)s'(t)-(y-s(t)))\frac{J_1(kr(w,t))}{r(w,t)},\\
L_2(w, t) &= L(w,t)-L_1(w,t)\ln{((x-t)^2+(y-s(t))^2)}.
\end{aligned}\end{equation}
Consequently, $M_1$, $M_2$, $L_1$, and $L_2$ are analytic functions. 

In this section, we focus on the single layer potential expressed as
\begin{equation}\begin{aligned}
S\psi(w) &=\int_{-1}^1 \ln{((x-t)^2+(y-s(t))^2)}\rho(w, t) dt+\int_{-1}^1 M_2(w,t)\psi(t) dt  \label{eq:single1},
\end{aligned}\end{equation}
where $\rho(w,t)=M_1(w,t) \psi(t)$. For each fixed target point $w$, $\rho(w,t)$ becomes a function purely of $t$. Thus,
we can apply the well-developed 1D Fourier extension technique~\cite{boyd2002comparison,bruno2007accurate,huybrechs2010fourier} 
to $\rho(w,t)$ to get
\begin{equation}
\rho(w,t) \approx \sum_{p=-P}^P a_pe^{\frac{\mathrm{i}p\pi}{2}t},
\end{equation}
where $a_p(w)$ can be computed by solving an optimization problem. Alternatively, the 2D Fourier extension 
\cite{matthysen2018function} algorithm can be applied directly to derive a 2D Fourier series approximation 
of $\rho(w,t)$ that is valid for all target points in the FMM leaf box. The optimization problem is usually exponentially 
ill-conditioned; however, the ill-conditioning can be resolved using the regularization technique described in 
\cite{barnett2022exponentially,huybrechs2010fourier}. Similar to the QB2X-Laplace technique \cite{ding2021quadrature}, 
we avoid the branch cut of the complex logarithmic function by using integration by parts to rewrite the single layer potential as
\begin{equation}\begin{aligned}
S\psi(w)&=\left[\ln{((x-t)^2+(y-s(t))^2)}f(t)\right]_{t=-1}^{t=1}\\
&+2\int_{-1}^1 \frac{(x-t)+(y-s(t))s'(t)}{(x-t)^2+(y-s(t))^2}f(t) dt+\int_{-1}^1 M_2(w,t)\psi(t) dt, \label{eq:single}
\end{aligned}\end{equation}
where $f(t)$ is the antiderivative of $\rho(w,t)$ given by
\begin{equation}
f(t) = \sum_{p=-P\\ p\neq 0}^P\frac{2a_p}{\mathrm{i}p\pi}e^{\frac{\mathrm{i}p\pi}{2}t}+a_0t
\end{equation}
for the given target point $w$. $f(t)$ can be further simplified by reapplying the Fourier extension
(or applying a precomputed mapping from polynomial basis to Fourier series) to the term $a_0 t$ to get
\begin{equation}\begin{aligned}
f(t) = \sum_{p=-P}^P b_p e^{\frac{\mathrm{i}p\pi}{2}t}.
\end{aligned}\end{equation}
In Eq. (\ref{eq:single}), the first term can be evaluated directly. The integrand in the last term is a smooth analytic 
function of $t$ and can be accurately and efficiently evaluated using any standard quadrature rules.
We focus on the numerical evaluation of the second term which contains the singular or near singular component of the kernel,
and rewrite the integral in the complex form as 
\begin{equation}\begin{aligned}
&2\int_{-1}^1 \frac{(x-t)f(t)+(y-s(t))s'(t)f(t)}{(x-t)^2+(y-s(t))^2} dt\\
=& \left(\int_{-1}^1 \frac{f(t)}{w-(t+is(t))}dt+\int_{-1}^1 \frac{f(t)}{w^*-(t-is(t))}dt \right)\\
&+ i\left(\int_{-1}^1 \frac{s'(t)f(t)}{w-(t+is(t))}dt-\int_{-1}^1 \frac{s'(t)f(t)}{w^*-(t-is(t))}dt\right),  \label{eq:complex}
\end{aligned}\end{equation}
where $w=x+\mathrm{i}y$ and $w^*=x-\mathrm{i}y$. We omit the details for the double layer potential as it can be simplified 
in an identical manner by replacing $M_1$ and $M_2$ by $L_1$ and $L_2$ in Eq. (\ref{eq:single1}), respectively.  Note that 
Eq.~(\ref{eq:complex}) is now similar to the layer potentials for the Laplace equation (e.g., see Eq. (6) in \cite{ding2021quadrature}).
Therefore, results from QB2X-Laplace can be adopted to derive the QB2X-Helmholtz representations.

\subsection{QB2X-Helmholtz for a straight line boundary}
To better understand the QB2X-Helmholtz technique, a straight line boundary is first considered ($s(t)=0$). 
The first two integrals in Eq. (\ref{eq:complex}) become complex integrals on the line segment from 
$-1$ to $1$ on the complex plane, both in the form of
\begin{equation}\begin{aligned}
\int_{-1}^1 \frac{f(z)}{w-z}\mathrm{d}z. \label{eq:simple_complex}
\end{aligned}\end{equation}
Since $s'(t)=0$, the last two integrals in Eq. (\ref{eq:complex}) become $0$.
We replace $f(z)$ in Eq.~(\ref{eq:simple_complex}) by its Fourier extension and separate the Fourier series 
into the non-negative and negative $p$ parts:
\begin{equation}\begin{aligned}
\int_{-1}^1 \frac{f(z)}{w-z}\mathrm{d}z&\approx\int_{-1}^1 \frac{\sum_{p=-P}^P b_p e^{\frac{\mathrm{i}p\pi}{2}z}}{w-z}\mathrm{d}z\\
&=\underbrace{\sum_{p=0}^P\int_{-1}^1 \frac{b_p e^{\frac{\mathrm{i}p\pi}{2}z} }{w-z}\mathrm{d}z}_{=I_1}+\underbrace{\sum_{p=-P}^{-1}\int_{-1}^1 \frac{ b_p e^{\frac{\mathrm{i}p\pi}{2}z}}{w-z}\mathrm{d}z}_{=I_2}.
\end{aligned}\end{equation}

\begin{figure}[bthp] 
   \centering
   \includegraphics[width=4.5in]{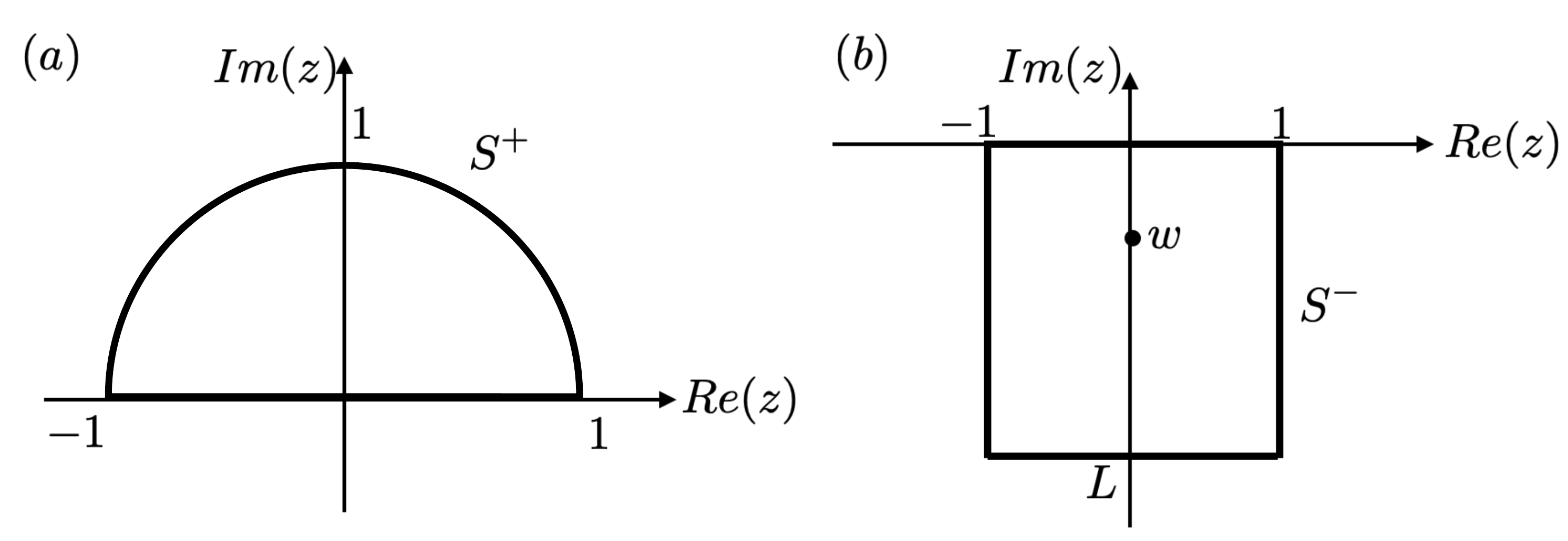} 
   \caption{Straight line boundary: (a) Contour $C_1$ for $I_1$ and (b) Contour $C_2$ for $I_2$}
   \label{fig:contourflat}
\end{figure}

The positive Fourier terms are bounded in the upper half-plane but grow exponentially in the lower half-plane 
when the imaginary part of $z \to -\infty$. For numerical stability considerations, two different 
contours are chosen to evaluate $I_1$ and $I_2$, respectively (See Fig.~\ref{fig:contourflat}).
We first use contour $C_1$ to compute $I_1$ which contains all the positive Fourier series terms. 
Because there is no pole inside the contour, the integral over $C_1$ is 0. As the integral over $C_1$ is the 
sum of the integral over $S^+$ and the line segment $[-1, 1]$, 
\begin{equation}\begin{aligned}
&\int_{C_1} \frac{b_p e^{\frac{\mathrm{i}p\pi}{2}z}}{w-z}\mathrm{d}z = \int_{S^+}\frac{b_p e^{\frac{\mathrm{i}p\pi}{2}z}}{w-z}\mathrm{d}z+\int_{-1}^{1}\frac{b_p e^{\frac{\mathrm{i}p\pi}{2}z}}{w-z}\mathrm{d}z  = 0\\
 &\Rightarrow \int_{-1}^1 \frac{b_p e^{\frac{\mathrm{i}p\pi}{2}z}}{w-z}\mathrm{d}z = -\int_{S^+} \frac{b_p e^{\frac{\mathrm{i}p\pi}{2}z}}{w-z}\mathrm{d}z.
\end{aligned}\end{equation}
Therefore, we derive a complex polynomial expansion for $I_1$ centered at $c$ as
\begin{equation}\begin{aligned}
&I_1 = \sum_{p=0}^{P}\int_{-1}^{1} \frac{b_p e^{\frac{\mathrm{i}p\pi}{2}z}}{w-z}\mathrm{d}z = -\sum_{p=0}^{P}\int_{S^+}\frac{b_p e^{\frac{\mathrm{i}p\pi}{2}z}}{w-z}\mathrm{d}z\\
&=-\sum_{p=0}^{P}\int_{S^+}\frac{b_p e^{\frac{\mathrm{i}p\pi}{2}z}}{(w-c)-(z-c)}\mathrm{d}z=\sum_{p=0}^{P}\int_{S^+}\frac{b_p e^{\frac{\mathrm{i}p\pi}{2}z}}{(z-c)} \frac{1}{(1-\frac{w-c}{z-c})}\mathrm{d}z\\
& = \sum_{p=0}^{P}\int_{S^+}\frac{b_p e^{\frac{\mathrm{i}p\pi}{2}z}}{(z-c)} \sum_{n=0}^\infty \left(\frac{w-c}{z-c}\right)^n \mathrm{d}z\\
&= \sum_{n=0}^\infty \left(\sum_{p=0}^{P}\int_{S^+}\frac{b_p e^{\frac{\mathrm{i}p\pi}{2}z}}{(z-c)^{n+1}}\mathrm{d}z\right)  (w-c)^n \approx  \sum_{n=0}^N u_n  (w-c)^n, \label{eq:I1_flat}
\end{aligned}\end{equation}
where
$u_n = \sum_{p=0}^{P}\int_{S^+}\frac{b_p e^{\frac{\mathrm{i}p\pi}{2}z}}{(z-c)^{n+1}}\mathrm{d}z$ are the local polynomial expansion coefficients for the upper half-plane.

Next, we use contour $C_2$ to compute $I_2$. Unlike the case for $C_1$, there exists a simple pole inside the counter $C_2$, namely, $w$. 
Therefore, the integral over $C_2$ can be computed using the Residue Theorem as follows:
\begin{equation}
\begin{aligned}
&\int_{C_2} \frac{b_p e^{\frac{\mathrm{i}p\pi}{2}z}}{w-z}\mathrm{d}z = \int_{S^-}\frac{b_p e^{\frac{\mathrm{i}p\pi}{2}z}}{w-z}\mathrm{d}z+\int_{-1}^{1}\frac{b_p e^{\frac{\mathrm{i}p\pi}{2}z}}{w-z}\mathrm{d}z \\ 
&=2\pi i Res[\frac{b_p e^{\frac{\mathrm{i}p\pi}{2}z}}{w-z}, w ] = 2\pi i b_pe^{\frac{\mathrm{i}p\pi}{2}w}\\
&\Rightarrow \int_{-1}^{1}\frac{b_p e^{\frac{\mathrm{i}p\pi}{2}z}}{w-z}\mathrm{d}z=-\int_{S^-}\frac{b_p e^{\frac{\mathrm{i}p\pi}{2}z}}{w-z}\mathrm{d}z+2\pi i b_pe^{\frac{\mathrm{i}p\pi}{2}w}.
\end{aligned}
\end{equation}
Applying the same separation of variables technique for the integral on $S^-$ as in Eq.~(\ref{eq:I1_flat}) for $I_1$, we have
\begin{equation}\begin{aligned}
&I_2= \sum_{p=-P}^{-1}\int_{-1}^1 \frac{b_p e^{\frac{\mathrm{i}p\pi}{2}z}}{w-z}\mathrm{d}z = -\sum_{p=-P}^{-1}\int_{S^{-}}\frac{b_p e^{\frac{\mathrm{i}p\pi}{2}z}}{w-z}\mathrm{d}z+2\pi i \sum_{p=-P}^{-1} b_p e^{\mathrm{i}pw}\\
&=\sum_{n=0}^\infty \left(\sum_{p=-P}^{-1}\int_{S^{-}}\frac{b_p e^{\frac{\mathrm{i}p\pi}{2}z}}{(z-c)^{n+1}}\mathrm{d}z\right)  (w-c)^n+2\pi i \sum_{p=-P}^{-1} b_p e^{\frac{\mathrm{i}p\pi}{2}w}\\
&\approx \sum_{n=0}^N d_n (w-c)^n+2\pi i \sum_{p=-P}^{-1} b_p e^{\frac{\mathrm{i}p\pi}{2}w},\label{eq:I2_flat}
\end{aligned}\end{equation}
where $d_n=\sum_{p=-P}^{-1}\int_{S^{-}}\frac{b_p e^{\frac{\mathrm{i}p\pi}{2}z}}{(z-c)^{n+1}}\mathrm{d}z$ are the local polynomial expansion coefficients for the lower half-plane.

In summary, by adding Eqs. (\ref{eq:I1_flat}) and (\ref{eq:I2_flat}), Eq. (\ref{eq:simple_complex}) 
can be approximated as the sum of a polynomial expansion and a plane wave type expansion
\begin{equation}\begin{aligned}
\int_{-1}^1 \frac{f(z)}{w-z}\mathrm{d}z&\approx\sum_{n=0}^N c_n  (w-c)^n+2\pi \mathrm{i} \sum_{p=-P}^{-1} b_p e^{\frac{\mathrm{i}p\pi}{2}w},
\end{aligned}\end{equation}
where 
\begin{equation}\begin{aligned}
c_n=\sum_{p=0}^{P}\int_{S^+}\frac{b_p e^{\frac{\mathrm{i}p\pi}{2}z}}{(z-c)^{n+1}}\mathrm{d}z+\sum_{p=-P}^{-1}\int_{S^{-}}\frac{b_p e^{\frac{\mathrm{i}p\pi}{2}z}}{(z-c)^{n+1}}\mathrm{d}z.
\end{aligned}\end{equation}
Since the source points on the segments $S^+$ and $S^-$ are well separated from the target point $w$, which is the only pole 
in the integrand, standard FMM error analysis can be applied to determine the necessary number of terms in the polynomial 
expansion. As shown in Sec.~\ref{sec:3}, for a prescribed accuracy requirement, the number of polynomial 
expansion terms needed is approximately the same as those from the worst-case FMM analysis for the 
Laplace kernel, i.e., $N=9$ for 3-digit, $N=18$ for 6-digit, $N=27$ for 9-digit, and $N=36$ for 12-digit accuracy.

\subsection{QB2X-Helmholtz for a curved boundary}
\begin{figure}[bthp] 
   \centering
   \includegraphics[width=4.5in]{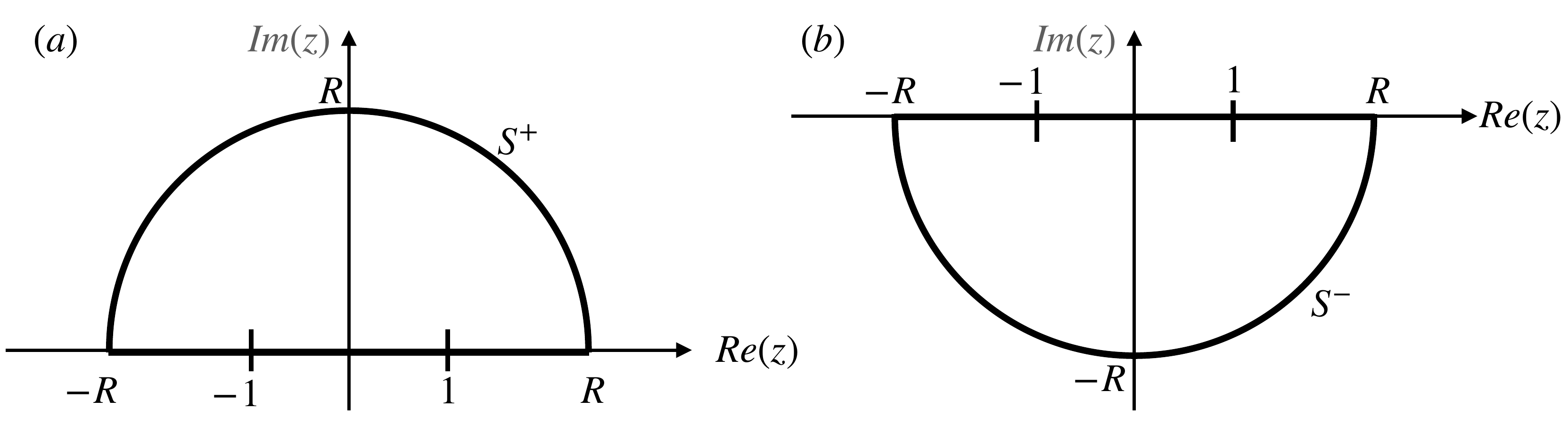} 
   \caption{Curved boundary: (a) $C_1$ contour for $I_1$ and (b) $C_2$ contour for $I_2$}
   \label{fig:contour}
\end{figure}

One main challenge for the QBX method is the slow convergence of the partial wave or polynomial expansions
when the boundary curve $s(t)$ is not flat or nearly-flat. This is demonstrated numerically in Sec.~\ref{sec:3}.
We present the QB2X-Helmholtz representation for a curved boundary which overcomes this challenge. We focus 
on the singular or nearly-singular parts of the Helmholtz layer potentials given by 
Eq.~(\ref{eq:complex}) in the complex form 
\begin{equation}\begin{aligned}
\int_{-1}^1 \frac{f(z)}{w-(z+\mathrm{i}s(z))}\mathrm{d}z.
	\label{eq:curved}
\end{aligned}\end{equation}
We follow the same strategy that was used for the straight line case. First, $f(z)$ is replaced by its Fourier extension 
and the Fourier series is separated into the non-negative and negative $p$ parts in order to ensure 
the convergence of the contour integrals over $C_1$ and $C_2$ (See Fig. \ref{fig:contour}), respectively.
\begin{equation}\begin{aligned}
\int_{-1}^1 \frac{f(z)}{w-(z+\mathrm{i}s(z))}\mathrm{d}z&=\int_{-1}^1 \frac{\sum_{p=-P}^P b_p e^{\frac{\mathrm{i}p\pi}{2}z}}{w-(z+\mathrm{i}s(z))}\mathrm{d}z\\
&\hskip-1in=\underbrace{\sum_{p=0}^P\int_{-1}^1 \frac{b_p e^{\frac{\mathrm{i}p\pi}{2}z} }{w-(z+\mathrm{i}s(z))}\mathrm{d}z}_{=I_1}+\underbrace{\sum_{p=-P}^{-1}\int_{-1}^1 \frac{ b_p e^{\frac{\mathrm{i}p\pi}{2}z}}{w-(z+\mathrm{i}s(z))}\mathrm{d}z}_{=I_2}.
\end{aligned}\end{equation}
Since the locations of the poles in the integrand of Eq.~(\ref{eq:curved}) are unknown, the contours $C_1$ and $C_2$ for 
the curved boundary are chosen differently from those for the straight line boundary. There are several advantages to the
current contour choice shown in Fig.~\ref{fig:contour}: Since $s(z)$ is a real-coefficient polynomial, there will be no pole on the 
real axis when $w$ is not located exactly on the boundary. Additionally, because of the polynomial growth in $s(z)$, all poles will 
be located inside either $C_1$ or $C_2$ if the radius $R$ in the contour is sufficiently large. Similar to Eq.~(16)
in \cite{ding2021quadrature}, this particular choice of contour allows for higher-order discretization of the boundary 
and easier control of the error when applying standard FMM error analysis for well-separated source and target points. 

The integral $I_1$ is computed using contour $C_1$ which is the union of four segments 
$S^+, [-1,1], [-R, -1]$, and $[1,R]$, see Fig.~\ref{fig:contour}. Applying the Residue Theorem,
\begin{equation}\begin{aligned}
&\int_{C_1} \frac{b_p e^{\frac{\mathrm{i}p\pi}{2}z}}{w-(z+\mathrm{i}s(z))}\mathrm{d}z = \left(\int_{S^+}+\int_{-1}^{1}+\int_{-R}^{-1}+\int_{1}^{R} \frac{b_p e^{\frac{\mathrm{i}p\pi}{2}z}}{w-(z+\mathrm{i}s(z))}\mathrm{d}z\right)\\
&\hspace{1cm} = 2\pi \mathrm{i} \sum_{j=1}^{J_u}Res[\frac{b_p e^{\frac{\mathrm{i}p\pi}{2}z}}{w-(z+\mathrm{i}s(z))}, r^u_j ]=-2\pi i \sum_{j=1}^{J_u} \frac{b_p e^{\frac{\mathrm{i}p\pi}{2}r^u_j}}{1+\mathrm{i}s'(r^u_j)}, \label{eq:residuethm1}
\end{aligned}\end{equation}
where $\{r^u_j\}_{j=1}^{J_u}$ are the roots of the polynomial equation $w-(z+\mathrm{i}s(z))=0$ in the upper 
complex half-plane. For a curved boundary $s(z) \neq 0$ with $s(0)=s'(0)=0$, when $R\rightarrow \infty$, 
the integral over $S^+$ approaches 0 for $p \geq 0$. Therefore, the integral $I_1$ can be computed as
\begin{equation}\begin{aligned}
I_1 &= \sum_{p=0}^{P}\int_{-1}^{1} \frac{b_p e^{\frac{\mathrm{i}p\pi}{2}z}}{w-(z+\mathrm{i}s(z))}\mathrm{d}z \\
&= -\sum_{p=0}^{P}\left( \int_1^\infty+\int_{-\infty}^{-1}\frac{b_p e^{\frac{\mathrm{i}p\pi}{2}z}}{w-(z+\mathrm{i}s(z))}\mathrm{d}z+2\pi \mathrm{i} \sum_{j=1}^{J_u} \frac{b_p e^{\frac{\mathrm{i}p\pi}{2}r^u_j}}{1+\mathrm{i}s'(r^u_j)}\right). \label{eq:splitradius1}
\end{aligned}\end{equation}
By applying the multipole/local expansion techniques from FMM~\cite{Carrier1988,greengard1987fast,greengard1997new}, 
the two integrals in $I_1$ from Eq.~(\ref{eq:splitradius1}) can be approximated using a complex Taylor expansion centered at $c$ as
\begin{equation}
\begin{aligned}
&-\sum_{p=0}^{P}\int_{1}^\infty+\int_{-\infty}^{-1}\frac{b_p e^{\frac{\mathrm{i}p\pi}{2}z}}{(w-c)-(z+\mathrm{i}s(z)-c)}\mathrm{d}z\\
&=\sum_{p=0}^{P}\int_1^\infty+\int_{-\infty}^{-1}\frac{b_p e^{\frac{\mathrm{i}p\pi}{2}z}}{(z+\mathrm{i}s(z)-c)(1-\frac{w-c}{(z+\mathrm{i}s(z)-c})}\mathrm{d}z\\
&= \sum_{p=0}^{P}\int_{1}^\infty+\int_{-\infty}^{-1}\frac{b_p e^{\frac{\mathrm{i}p\pi}{2}z}}{(z+\mathrm{i}s(z)-c)} \sum_{n=0}^\infty \left(\frac{w-c}{z+\mathrm{i}s(z)-c}\right)^n \mathrm{d}z\\
&= \sum_{n=0}^\infty \left(\sum_{p=0}^{P}\int_{1}^{\infty}+\int_{-\infty}^{-1}\frac{b_p e^{\frac{\mathrm{i}p\pi}{2}z}}{(z+\mathrm{i}s(z)-c)^{n+1}}\mathrm{d}z\right)  (w-c)^n\\
& \approx  \sum_{n=0}^N u_n  (w-c)^n,
\end{aligned}
\end{equation}
where the number of terms in the expansion is determined by standard FMM error analysis.
Therefore, the QB2X representation of the integral $I_1$ is given by
\begin{equation}\begin{aligned}
I_1\approx  \sum_{n=0}^N u_n  (w-c)^n-2\pi \mathrm{i}\sum_{p=0}^P  \sum_{j=1}^{J_u} \frac{b_p e^{\frac{\mathrm{i}p\pi}{2}r^u_j}}{1+\mathrm{i}s'(r^u_j)},
\end{aligned}\end{equation}
where the polynomial expansion coefficients are
\begin{equation}\begin{aligned}
u_n=\sum_{p=0}^{P}\left(\int_{1}^{\infty}+\int_{-\infty}^{-1}\frac{b_p e^{\frac{\mathrm{i}p\pi}{2}z}}{(z+\mathrm{i}s(z)-c)^{n+1}}\mathrm{d}z\right).
\end{aligned}\end{equation}

Similarly, we use contour $C_2=S^- \cup [-1,1] \cup [-R, -1] \cup [1,R]$ to compute the integral $I_2$. Using the Residue Theorem
\begin{equation}\begin{aligned}
&\int_{C_2} \frac{b_p e^{\frac{\mathrm{i}p\pi}{2}z}}{w-(z+\mathrm{i}s(z))}\mathrm{d}z = (\int_{S^-}+\int_{-1}^{1}+\int_{-R}^{-1}+\int_{1}^{R} \frac{b_p e^{\frac{\mathrm{i}p\pi}{2}z}}{w-(z+\mathrm{i}s(z))}\mathrm{d}z)\\
&~~~~~~~ = -2\pi \mathrm{i} \sum_{j=1}^{J_d}Res[\frac{b_p e^{\frac{\mathrm{i}p\pi}{2}z}}{w-(z+\mathrm{i}s(z))}, r^d_j ]=2\pi i \sum_{j=1}^{J_d} \frac{b_p e^{\frac{\mathrm{i}p\pi}{2}r^d_j}}{1+\mathrm{i}s'(r^d_j)},\label{eq:residuethm2}
\end{aligned}\end{equation}
where $\{r^d_j\}_{j=1}^{J_d}$ are the roots of the polynomial equation $w-(z+\mathrm{i}s(z))=0$ in the lower complex half-plane. 
When $R\rightarrow \infty$, the integral over $S^-$ approaches 0, so the integral $I_2$ can be computed by
\begin{equation}\begin{aligned}
I_2 &= \sum_{p=-P}^{-1}\int_{-1}^{1} \frac{b_p e^{\frac{\mathrm{i}p\pi}{2}z}}{w-(z+\mathrm{i}s(z))}\mathrm{d}z \\
&= -\sum_{p=-P}^{-1}\left( \int_1^\infty+\int_{-\infty}^{-1}\frac{b_p e^{\frac{\mathrm{i}p\pi}{2}z}}{w-(z+\mathrm{i}s(z))}\mathrm{d}z-2\pi \mathrm{i} \sum_{j=1}^{J_d} \frac{b_p e^{\frac{\mathrm{i}p\pi}{2}r^d_j}}{1+\mathrm{i} s'(r^d_j)}\right). \label{eq:splitradius2}
\end{aligned}\end{equation}
Applying the same FMM technique, the two integrals in $I_2$ from Eq.~(\ref{eq:splitradius2}) can be approximated 
using a complex Taylor expansion centered at $c$ as
\begin{equation}\begin{aligned}
&-\sum_{p=-P}^{-1}\int_{1}^\infty+\int_{-\infty}^{-1}\frac{b_p e^{\frac{\mathrm{i}p\pi}{2}z}}{(w-c)-(z+\mathrm{i}s(z)-c)}\mathrm{d}z\\
&~~~~~~= \sum_{n=0}^\infty \left(\sum_{p=-P}^{-1}\int_{1}^{\infty}+\int_{-\infty}^{-1}\frac{b_p e^{\frac{\mathrm{i}p\pi}{2}z}}{(z+\mathrm{i}s(z)-c)^{n+1}}\mathrm{d}z\right)  (w-c)^n\\
&~~~~~~ \approx  \sum_{n=0}^N d_n  (w-c)^n.
\end{aligned}\end{equation}
Therefore, the QB2X representation of $I_2$ is given by
\begin{equation}\begin{aligned}
I_2\approx  \sum_{n=0}^N d_n  (w-c)^n+2\pi \mathrm{i}\sum_{p=-P}^{-1}  \sum_{j=1}^{J_d} \frac{b_p e^{\frac{\mathrm{i}p\pi}{2}r^d_j}}{1+\mathrm{i}s'(r^d_j)},
\end{aligned}\end{equation}
where the polynomial expansion coefficients are
\begin{equation}\begin{aligned}
d_n=\sum_{p=-P}^{-1}\left(\int_{1}^{\infty}+\int_{-\infty}^{-1}\frac{b_p e^{\frac{\mathrm{i}p\pi}{2}z}}{(z+\mathrm{i}s(z)-c)^{n+1}}\mathrm{d}z\right).
\end{aligned}\end{equation}

We summarize our results by combining $I_1$ and $I_2$ and present the QB2X representation of Eq.~(\ref{eq:curved}) as
\begin{equation}\begin{aligned}\label{eq:QB2X}
\int_{-1}^1 \frac{f(z)}{w-(z+\mathrm{i}s(z))}\mathrm{d}z & \\ 
&\hskip-1.5in\approx \sum_{n=0}^N c_n (w-c)^n-2\pi \mathrm{i} \left( \sum_{p=0}^P  \sum_{j=1}^{J_u} \frac{b_p e^{\frac{\mathrm{i}p\pi}{2}r^u_j}}{1+\mathrm{i}s'(r^u_j)}-\sum_{p=-P}^{-1}  \sum_{j=1}^{J_d} \frac{b_p e^{\frac{\mathrm{i}p\pi}{2}r^d_j}}{1+\mathrm{i}s'(r^d_j)}\right),
\end{aligned}\end{equation} 
where 
\begin{equation}\begin{aligned}
&c_n = u_n+d_n\\
&=\sum_{p=-P}^{P}b_p\left(\int_{1}^{\infty}\frac{e^{\frac{\mathrm{i}p\pi}{2}z}}{(z+\mathrm{i}s(z)-c)^{n+1}}\mathrm{d}z+\int_{-\infty}^{-1}\frac{e^{\frac{\mathrm{i}p\pi}{2}z}}{(z+\mathrm{i}s(z)-c)^{n+1}}\mathrm{d}z\right).
\end{aligned}\end{equation}
The singular or nearly-singular part of the Helmholtz layer potential in Eq. (\ref{eq:complex}) is simply a combination of 
four terms of Eq. (\ref{eq:QB2X}) with different density functions. Note that the integrals composing $c_n$ are independent of 
the target points. For better efficiency, these integrals are only computed once for each FMM tree leaf and stored in 
a two dimensional table for different $p$ and $n$ values. This table is then used for all target points.

\subsection{Stable computation of the QB2X-Helmholtz representation} 
In the numerical computation of the QB2X-Helmholtz representation, when some of the roots of the polynomial equation 
$w-(z+\mathrm{i}s(z))=0$ are close to each other in rare cases, the computation of the residues (plane wave type expansion) in 
Eqs. (\ref{eq:residuethm1}) and (\ref{eq:residuethm2})
\begin{equation}\begin{aligned}
&\frac{1}{2\pi i}\int_{C_1} \frac{b_p e^{\frac{\mathrm{i}p\pi}{2}z}}{w-(z+\mathrm{i}s(z))}\mathrm{d}z \mbox{ and }\frac{1}{2\pi i}\int_{C_2} \frac{b_p e^{\frac{\mathrm{i}p\pi}{2}z}}{w-(z+\mathrm{i}s(z))}\mathrm{d}z \label{eq:stabledirect}
\end{aligned}\end{equation}
may experience a loss of significance due to the subtraction of two close numbers. For example, consider a third degree 
polynomial $s(z)$ with leading coefficient $c_0$, and assume the denominator $(w-(z+\mathrm{i}s(z)))=0$ has three roots: $r_1$, $r_2$ in the upper complex half-plane that are very close to each other, and $r_3$ in the lower complex half-plane. Then the plane wave type expansion from the Residue Theorem becomes
\begin{equation}\begin{aligned}
&\frac{1}{2\pi i}\int_{C_1} \frac{b_p e^{\frac{\mathrm{i}p\pi}{2}z}}{c_0 (z-r_1)(z-r_2)(z-r_3)}\mathrm{d}z\\
&= \frac{b_p e^{\frac{\mathrm{i}p\pi}{2}r_1}}{c_0(r_1-r_2)(r_1-r_3)}+\frac{b_p e^{\frac{\mathrm{i}p\pi}{2}r_2}}{c_0(r_2-r_1)(r_2-r_3)}, \label{eq:directResidue}
\end{aligned}\end{equation}
where $C_1$ is the contour given in Fig. \ref{fig:contour}(a). Eq.~(\ref{eq:directResidue}) clearly shows the loss of accuracy when roots are close to each other. We present a remedy  for this issue. Let $r_c$ be the midpoint of $r_1$ and $r_2$ and factor out $z-r_c$ from the denominator as
\begin{align}
&\int_{C_1} \frac{b_p e^{\frac{\mathrm{i}p\pi}{2}z}}{c_0 (z-r_1)(z-r_2)(z-r_3)}\mathrm{d}z= \int_{C_1} \frac{b_p e^{\frac{\mathrm{i}p\pi}{2}z}}{c_0 (z-r_c)^2 \left(\frac{z-r_1}{z-r_c}\right)\left(\frac{z-r_2}{z-r_c}\right) (z-r_3)} \mathrm{d}z \nonumber\\
&=\int_{C_1} \frac{b_p e^{\frac{\mathrm{i}p\pi}{2}z}}{c_0 (z-r_c)^2 \left(\frac{z-r_c+r_c-r_1}{z-r_c}\right)\left(\frac{z-r_c+r_c-r_2}{z-r_c}\right)(z-r_3)} \mathrm{d}z \nonumber\\
&=\int_{C_1} \frac{b_p e^{\frac{\mathrm{i}p\pi}{2}z}}{c_0 (z-r_c)^2 \left(1-\frac{r_1-r_c}{z-r_c}\right)\left(1-\frac{r_2-r_c}{z-r_c}\right)(z-r_3)} \mathrm{d}z.
\end{align}
Using the geometric series of $1/\left(1-\frac{r_1-r_c}{z-r_c}\right)$ and $1/\left(1-\frac{r_2-r_c}{z-r_c}\right)$:
\begin{align}
&=\int_{C_1} \frac{b_p e^{\frac{\mathrm{i}p\pi}{2}z}}{c_0 (z-r_c)^2 (z-r_3)}\left(1+\frac{r_1-r_c}{z-r_c}+\left(\frac{r_1-r_c}{z-r_c}\right)^2+\cdots\right)\nonumber\\
&~~~~~~~~~~~~~~~~~~~~~~~~~~~~~~~~~~\times \left(1+\frac{r_2-r_c}{z-r_c}+\left(\frac{r_2-r_c}{z-r_c}\right)^2+\cdots\right) \mathrm{d}z\nonumber\\
&=\int_{C_1} \frac{b_p e^{\frac{\mathrm{i}p\pi}{2}z}}{c_0 (z-r_c)^2(z-r_3) }\left(\delta_0+\frac{\delta_1}{z-r_c}+\frac{\delta_2}{(z-r_c)^2}+\cdots\right) \mathrm{d}z\label{eq:deltaCoefs}\\
&=\frac{b_p}{c_0}\int_{C_1} \frac{\delta_0e^{\frac{\mathrm{i}p\pi}{2}z}}{(z-r_c)^2 (z-r_3)}+\frac{\delta_1 e^{\frac{\mathrm{i}p\pi}{2}z}}{(z-r_c)^3(z-r_3) }+\frac{\delta_2  e^{\frac{\mathrm{i}p\pi}{2}z}}{(z-r_c)^4 (z-r_3)}+\cdots \mathrm{d}z, \nonumber
\end{align}
where $\{\delta_j\}_{j=0}^{\infty}$ are the coefficients of the product of the two geometric series. 
By applying the Residue Theorem for each term and truncating the sum at $M$th term, the integral can be approximated by 
\begin{align}
\int_{C_1} \frac{b_p e^{\frac{\mathrm{i}p\pi}{2}z}}{c_0 (z-r_1)(z-r_2)(z-r_3)}\mathrm{d}z \approx\frac{2 \pi i b_p}{c_0}  \sum_{m=1}^{M} \delta_{m-1} \left.\left(\frac{e^{\frac{ip\pi}{2}z}}{z-r_3}\right)^{(m)}\right|_{z=r_c}. \label{eq:accurate}
\end{align}
Since the geometric series decays rapidly, only a small number of terms is required in this alternative representation 
for cases when roots are close to each other.

We present comparisons of the original plane wave type formula in Eq.~(\ref{eq:directResidue}) with the alternative representation in
Eq.~(\ref{eq:accurate}) for a curved boundary given by $s(z)=-\left(\frac{4\times 10^{19}}{-8.7\times 10^{19}+\pi^2}\right)z^3+\left(\frac{-6.3\times 10^{19}+\pi^2}{-8.7\times 10^{19}+\pi^2}\right)z^2  $ and $w=\left(\frac{3}{5}+\frac{29}{100}\mathrm{i}\right)+\frac{\mathrm{i} \pi^2}{10^{20}}+\frac{4.48\times 10^{19}+6.4\times 10^{19} \mathrm{i}}{-8.7\times 10^{19}+\pi^2}$.
Three roots of $w-(z+\mathrm{i}s(z))=0$ are $r_1=\frac{1}{5}-\pi\times10^{-10}+\frac{1}{2}\mathrm{i}$, $r_2=\frac{1}{5}+\pi\times10^{-10}+\frac{1}{2}\mathrm{i}$, and $r_3=\frac{-7.8\times 10^{19}+\pi^2}{4\times 10^{19}}-\mathrm{i}$ and construct the Helmholtz 
layer potential using $p=1$ and $b_p=1$.  
Integration over $C_1$ using the Residue Theorem in Eq. (\ref{eq:directResidue}) yields a result of 
$0.6757127940654755-0.07931315898895264\mathrm{i}$ and the new representation in Eq. (\ref{eq:accurate}) with $M=3$ yields 
$0.67571272714158517293 - 0.07931328249252365183\mathrm{i}$. Eq. (\ref{eq:stabledirect}) is evaluated using {\tt NIntegrate} 
command with a precision goal of 24 digits in Mathematica  as the reference solution. The absolute error between the reference solution and the direct Residue Theorem and the new representation is $1.4\times 10^{-7}$ and is $4.6\times 10^{-18}$, respectively. An additional example with a fourth degree polynomial can be found in \cite{ding2021quadrature}. Note that loss of accuracy does not happen in the case of second degree polynomials due to the symmetry of the two roots.

\section{Numerical Results}
\label{sec:3}
In this section, we present numerical results to demonstrate the effectiveness of the QB2X technique for the Helmholtz
layer potentials. We test the accuracy of QB2X-Helmholtz for different $P$, $N$, and $k$ values, where $P$ is the 
number of terms used for the Fourier extensions, $N$ is the number of terms in the complex Taylor expansions, 
and $k$ is the wave number. In this paper, we restrict our attention only to the low- to mild-frequency $k$ values 
which is required when evaluating the local direct interactions in the FMM algorithm. We compute
reference solutions using the {\tt NIntegrate} command from Mathematica with a precision goal of 24 digits for 
all numerical examples. We use the $L^\infty$ norm when comparing QB2X or QBX solutions to the reference solutions, 
and present the number of digits in the error analysis computed using 
$\|\log|err|\|_{\infty} = \max\limits_{w}\log\|Q\psi(w) - R\psi(w)\|$, where $Q\psi(w)$ is the QBX or QB2X solution 
and $R\psi(w)$ is the reference Helmholtz single layer or double layer solution.
 \begin{figure}[htbp] 
   \centering
   \includegraphics[width=4.5in]{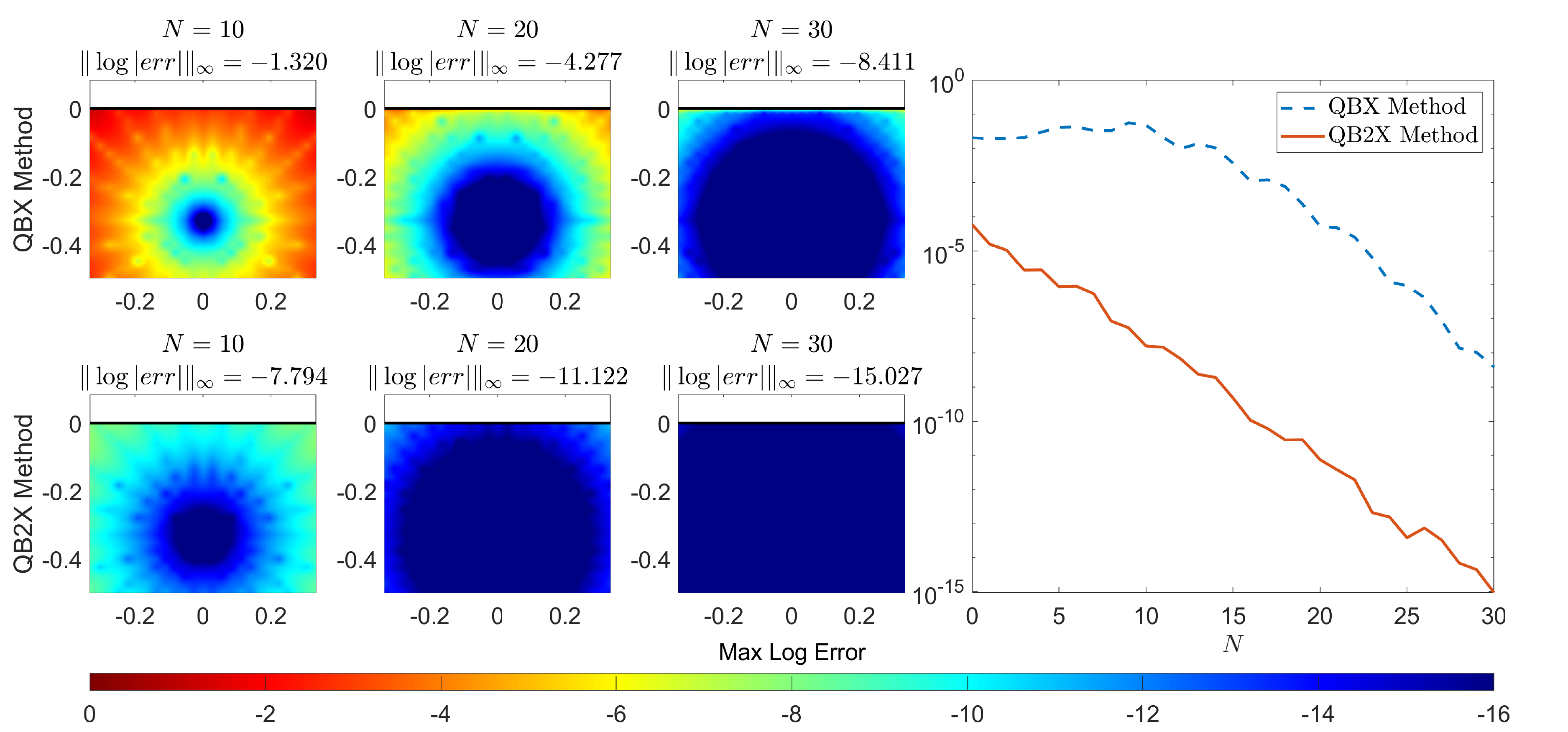} 
	 \caption{Single layer potentials for $s(t)=0$: (\textit{Left}) Errors from both QBX (upper) and QB2X (lower)
	 with respect to reference solutions for density function $\psi(t)=\cos(20t)$ with $P=50$. 
	 (\textit{Right}) Convergence of the $L^{\infty}$ errors for both the QBX (dashed line) 
	 and QB2X (solid line) cases with respect to $N$.}
   \label{fig:SLP_s0_QBX_QB2X}
\end{figure}

In the first example, we compare the original QBX-Helmholtz method with the new QB2X-Helmholtz method for a flat 
boundary $s(t)=0$ with a mildly oscillating trigonometric density function $\psi(t)=\cos(20t)$. We set $P=50$ and $k=1$.
Fig. \ref{fig:SLP_s0_QBX_QB2X} shows the errors of both the QBX and QB2X methods for $N=10, 20$ and $30$ and the 
convergence for different $N$ values. For this example, QB2X-Helmholtz reaches its asymptotic error minimum when 
$N=31$, and $\|\log|err|\|_{\infty} = -15.033$ which is approximately machine precision.
\begin{figure}[htbp] 
   \centering
   \includegraphics[width=4.0in]{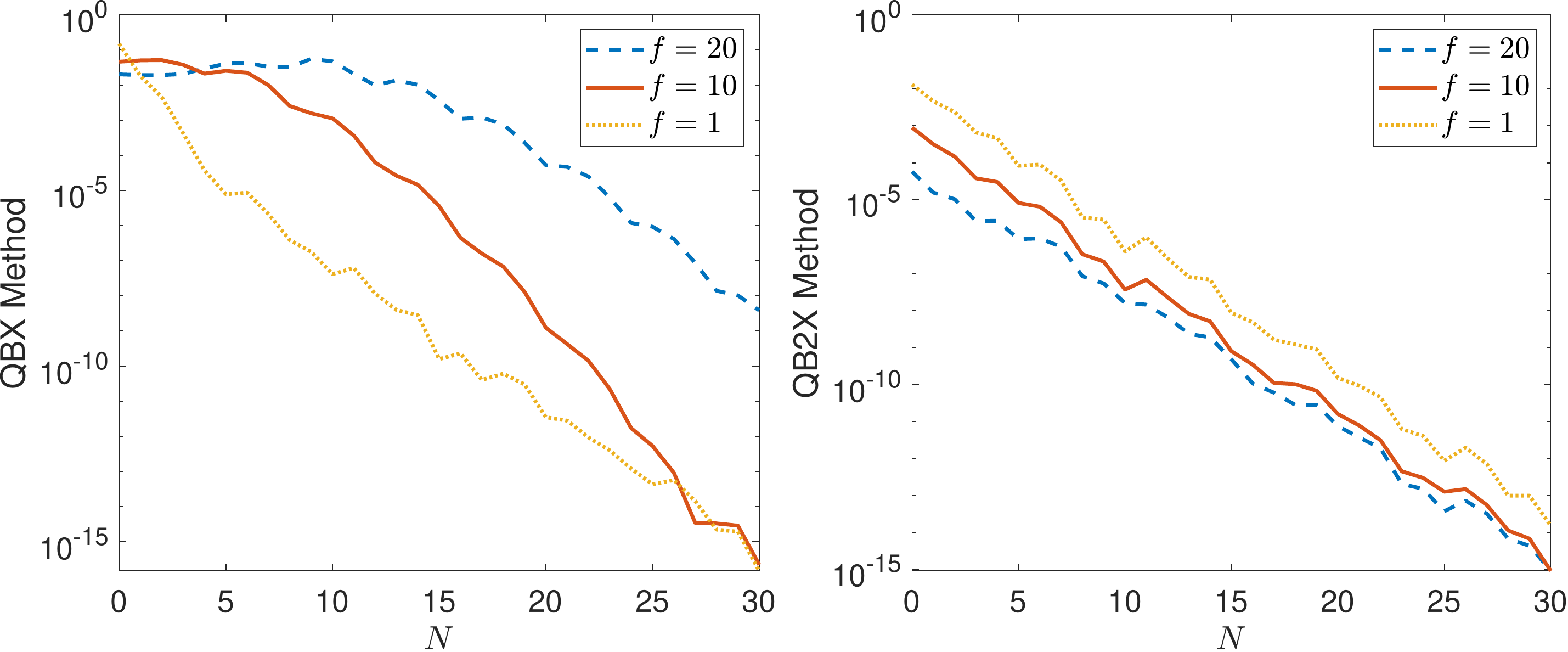} 
   \caption{Single layer potentials when $s(t)=0$: (\textit{Left}) QBX and (\textit{Right}) QB2X convergence for 
	density function $\psi(t)=\cos(ft)$ with $f$=1 (dotted line), 10 (solid line), 20 (dashed line) and 
	increasing $N$ values. We set $P = 50$.}
   \label{fig:SLP_s0_rho_phase}
\end{figure}

When the boundary is a straight line segment, we noticed that the convergence of QBX-Helmholtz depends on the periodicity (frequency) of the density
function. For low-frequency density functions, both QBX-Helmholtz and QB2X-Helmholtz work well. However as the frequency of the density function
increases, the performance of QBX deteriorates rapidly. On the other hand, the error in QB2X follows the standard FMM error analysis and the
convergence (with respect to $N$) is almost independent of the frequency of the density function. This is demonstrated in 
Figure~\ref{fig:SLP_s0_rho_phase}  for different frequencies $f$ in the density function $\psi(t) = \cos(ft)$.
\begin{figure}[htbp] 
   \centering
   \includegraphics[width=4.5in]{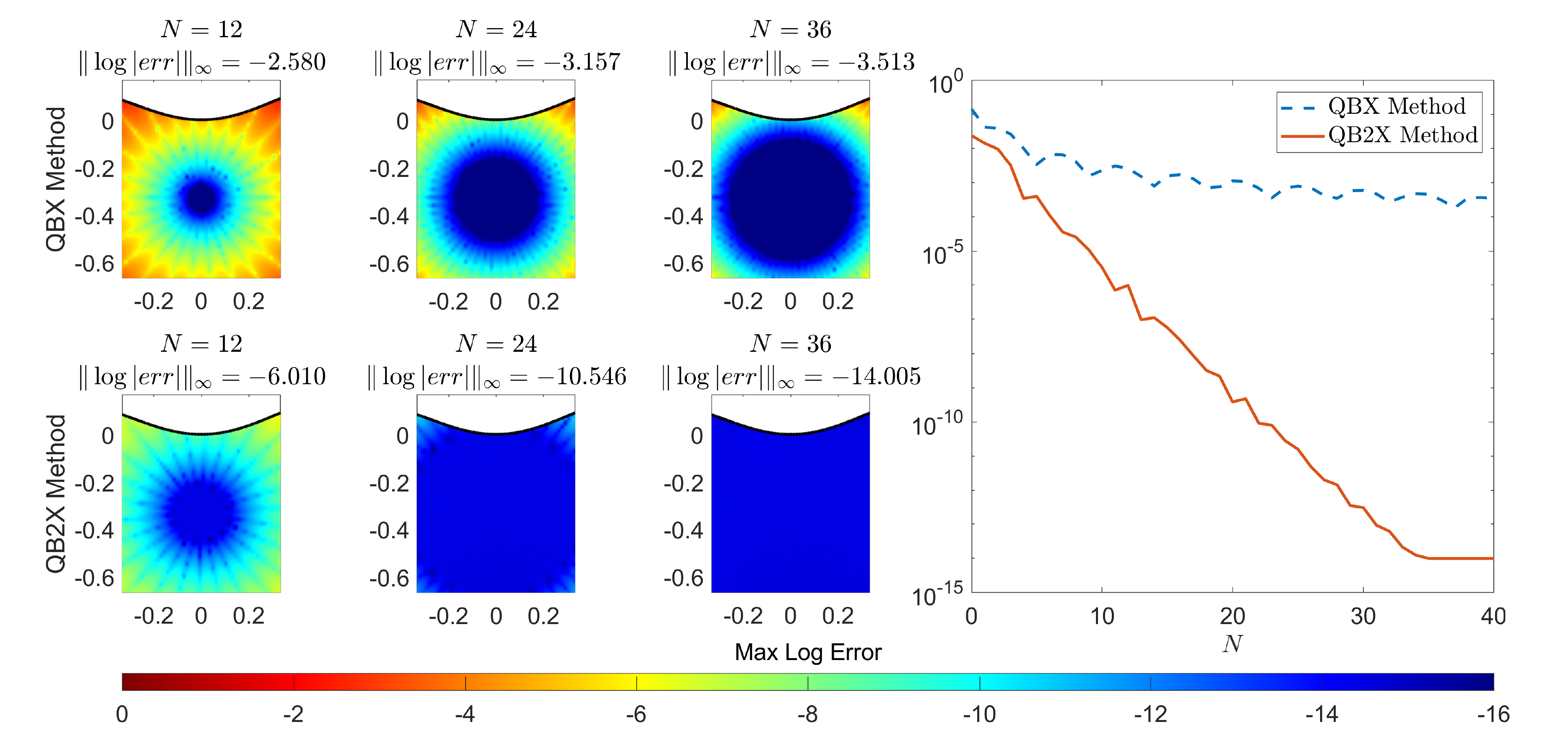} 
	\caption{Single layer potentials for $s(t)=t^2+\frac{t^3}{10}-2t^4$: (\textit{Left}) Errors from QBX (upper) 
	and QB2X (lower) for density function $\psi(t)=1$ with $P=100$. (\textit{Right}) Convergence of the 
	$L^{\infty}$ error in QBX (dashed line) and QB2X (solid line) for increasing $N$ values.}
   \label{fig:SLP_smcurve_QBX_QB2X}
\end{figure}

In the second example, we consider a non-flat boundary given by $s(t)=t^2+\frac{t^3}{10}-2t^4$. 
We focus on the nonlinear contribution from the boundary geometry and set the density function to the 
constant function $\psi(t)=1$ and $k=1$. In the left of Figure~\ref{fig:SLP_smcurve_QBX_QB2X}, we 
show the error for $N=12$, $24$, and $36$. For the QBX representation centered at $(0, -\frac{1}{3})$, 
although the $L^\infty$ error is small close to the center, the error close to the boundary and 
FMM box edges decays very slowly when $N$ increases. On the other hand, the QB2X method achieves machine 
precision accuracy at all points in the FMM leaf box when $N=36$. In the right of Figure~\ref{fig:SLP_smcurve_QBX_QB2X},
we plot the convergence of QBX and QB2X as a function of $N$. The convergence of QB2X follows the standard FMM 
error analysis and reaches machine precision when $N \approx 36$. The QBX method is unable to effectively resolve 
the nonlinearity from the curved boundary and the error decays slowly. In both cases, $P=100$ is used in the 
Fourier extension.

\begin{figure}[htbp] 
   \centering
   \includegraphics[width=4.5in]{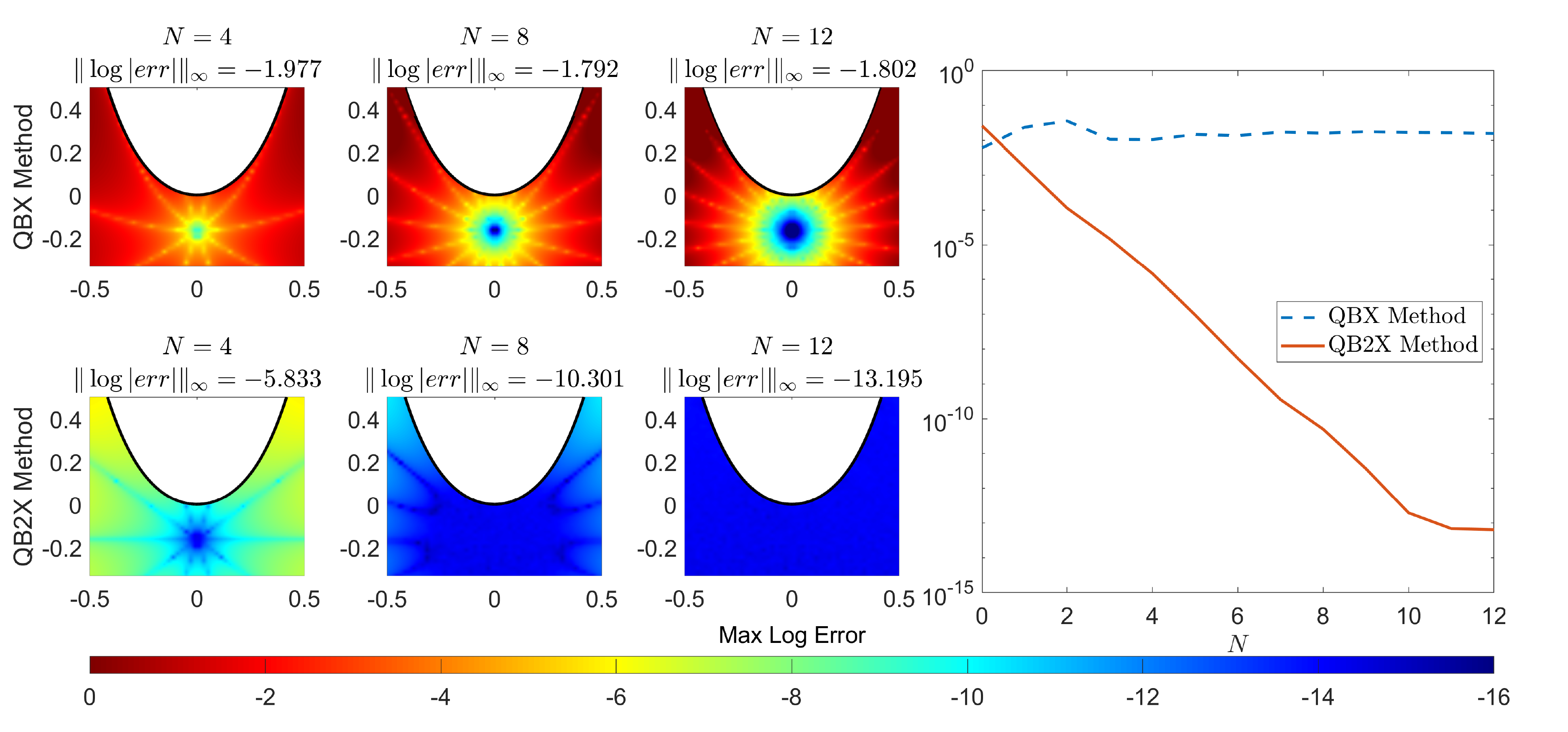} 
   \caption{Single layer potentials when $s(t)=5t^4+2t^2$: (\textit{Left}) Errors from QBX and QB2X for $\psi(t)=1$, $N=4$, $8$, $12$, and $P=100$. 
	(\textit{Right}) $L^{\infty}$ error convergence for both QBX (dashed line) and QB2X (solid line) for increasing $N$ values.}
   \label{fig:SLP_parabola_QBX_QB2X}
\end{figure}

In order to further demonstrate the convergence of the QBX and QB2X methods and its dependency on the boundary geometry, 
in the third example, we consider the boundary setting $s(t)=5t^4+2t^2$ which curves rapidly away from the partial 
wave (QBX) or polynomial (QB2X) expansion center. This extreme geometry setting should never happen in a uniform FMM tree structure 
but may appear when an adaptive FMM tree is used to more effectively resolve the solution. 
All the other features are the same as those in the previous example. 
On the left of Fig. \ref{fig:SLP_parabola_QBX_QB2X}, the error from QBX decays slowly in most of the computational
domain. This is not surprising as most existing QBX implementations require the resolved boundary in a leaf box
of the FMM tree is ``close" to flat by oversampling the boundary geometry. On the other hard, the QB2X method allows
high degree polynomial description of the boundary, and the error in this example converges to approximately machine 
precision when $N=11$. On the right of Fig.~\ref{fig:SLP_parabola_QBX_QB2X}, we present the $L^\infty$ error convergence
for both methods. Clearly, the convergence of QB2X follows the standard FMM error analysis. For the particular geometry setting, as the FMM leaf box size is $\frac{2}{3} \times \frac{2}{3}$, we have $\|w-c\| \leq \frac{\sqrt{2}}{3}$, and for $z \in [-\infty, -1] \cup [1, \infty]$, $\| z +is(z)– c \| \geq \sqrt{(\pm 1-0)^2 +(s(\pm 1)-0)^2}$. Therefore the local polynomial expansion converges like $\left( \frac{\|w-c\|}{\| z+is(z)–c \| } \right)^N \approx (1/15)^N$. Applying this worst case analysis, when $N \approx 12$, $(1/15)^{12} \approx 7.71\times 10^{-15}$, QB2X should achieve machine precision accuracy. This analysis is consistent with the numerical results.

\begin{figure}[htbp] 
   \centering
   \includegraphics[width=4.5in]{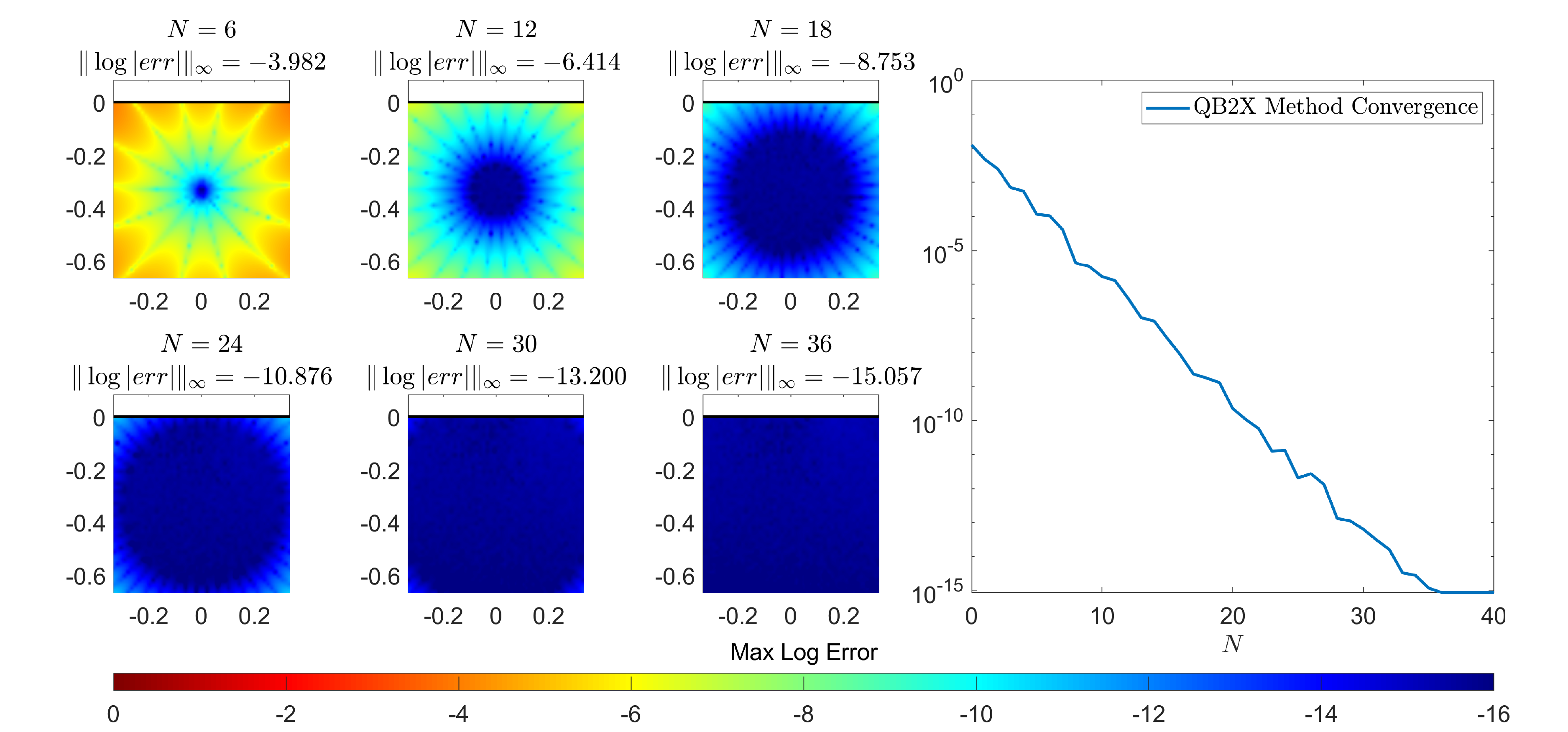} 
   \caption{Single layer potentials for $s(t)=0$ and $\psi(t)=\frac{2t^2+2t+3}{4}$. $P=50$: (\textit{Left}) QB2X errors for 
	$N=6$, $12$, $18$, $24$, $30$, and $36$. $P=50$. (\textit{Right}) Convergence of the $L^\infty$ error for different $N$ values.}
   \label{fig:SLP_s0_rhopoly}
\end{figure}

In many existing implementations of the boundary integral equation method, the density functions are usually 
approximated by piecewise polynomials (instead of a Fourier series). We therefore consider a polynomial 
density function $\psi(t)=\frac{2t^2+2t+3}{4}$ in the next example and show the performance of QB2X. 
In Figure \ref{fig:SLP_s0_rhopoly}, we show the errors for $N=6$, $12$, $18$, $24$, $30$, and $36$ and 
convergence of the QB2X method for a flat boundary. We set $P=50$ and $k=1$. The numerical results show that 
QB2X reaches machine precision accuracy at about $N=36$ where $\|\log| err |\|_{\infty} = -15.0575$.

There are at least two ways to derive the Fourier series required in the QB2X method, by applying existing Fourier
extension algorithms on the fly, or first representing the density functions as piecewise polynomials, and then
map the polynomials to Fourier series using precomputed (using Fourier extension) transformation tables from 
the polynomial basis to Fourier basis. The latter technique can be more efficient but requires additional storage. 

\begin{figure}[htbp] 
   \centering
   \includegraphics[width=4.5in]{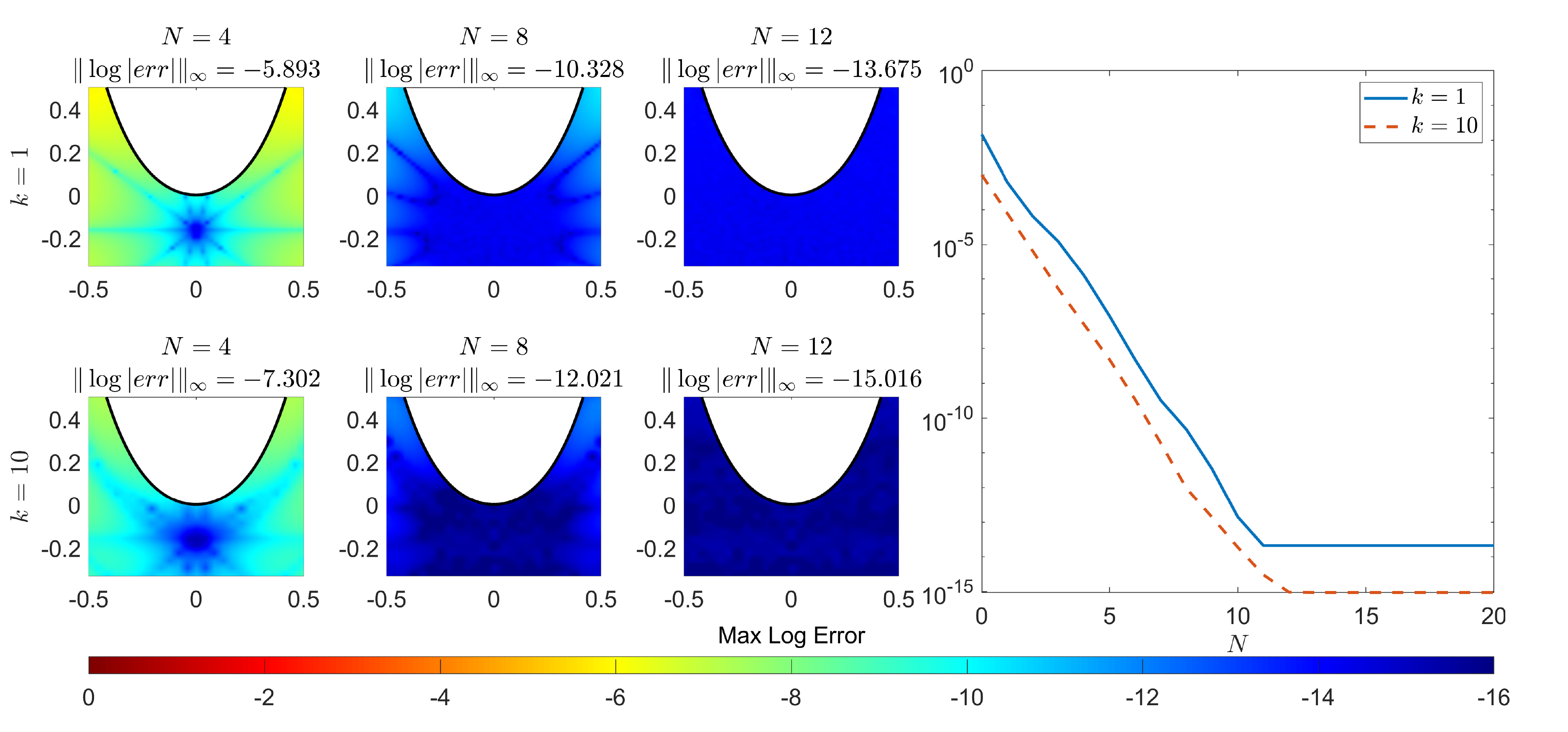} 
   \caption{Single layer potentials for $s(t) = 2t^2+5t^4$ and $\psi(t)=1$: (\textit{Left}) QB2X errors for wave numbers $k=1$ and $P=100$ (upper),
	and for $k=10$ and $P=300$ (lower) when $N=4$, $8$, and $12$. 
	(\textit{Right}) QB2X convergence ($L^\infty$ error) as a function of $N$ for 
	$k=1$ (solid line) and $k=10$ (dashed line).}
   \label{fig:SLP_parabola_kvary}
\end{figure}

In the FMM algorithm, the evaluation of the singular and nearly-singular Helmholtz layer potentials is only required 
in the local direct interactions of the leaf boxes of the FMM tree structure. Therefore, we only need to consider 
low- and mild-frequency $k$ values. In the next example, we consider the curved boundary $s(t)=5t^4+2t^2$ and 
density function $\psi(t)=1$ for wave numbers $k=1$ and $k=10$. Clearly, for larger $k$ value,
more Fourier extension terms are required to capture the linear and nonlinear contributions. 
We set $P=100$ for $k=1$ and $P=300$ for $k=10$. On the left of Fig. \ref{fig:SLP_parabola_kvary}, we show 
the error results for $N=4$, $8$, and $12$ and different $k$ values. On the right, we show the convergence of the
QB2X-Helmholtz method as $N$ increases. The numerical results show that when a sufficient number of Fourier series terms 
are used, the convergence properties for $k=1$ and $k=10$ as a function of $N$ are at similar rates, and reach 
machine precision for all target points inside the FMM leaf box. 

\begin{figure}[htbp] 
   \centering
   \includegraphics[width=4.5in]{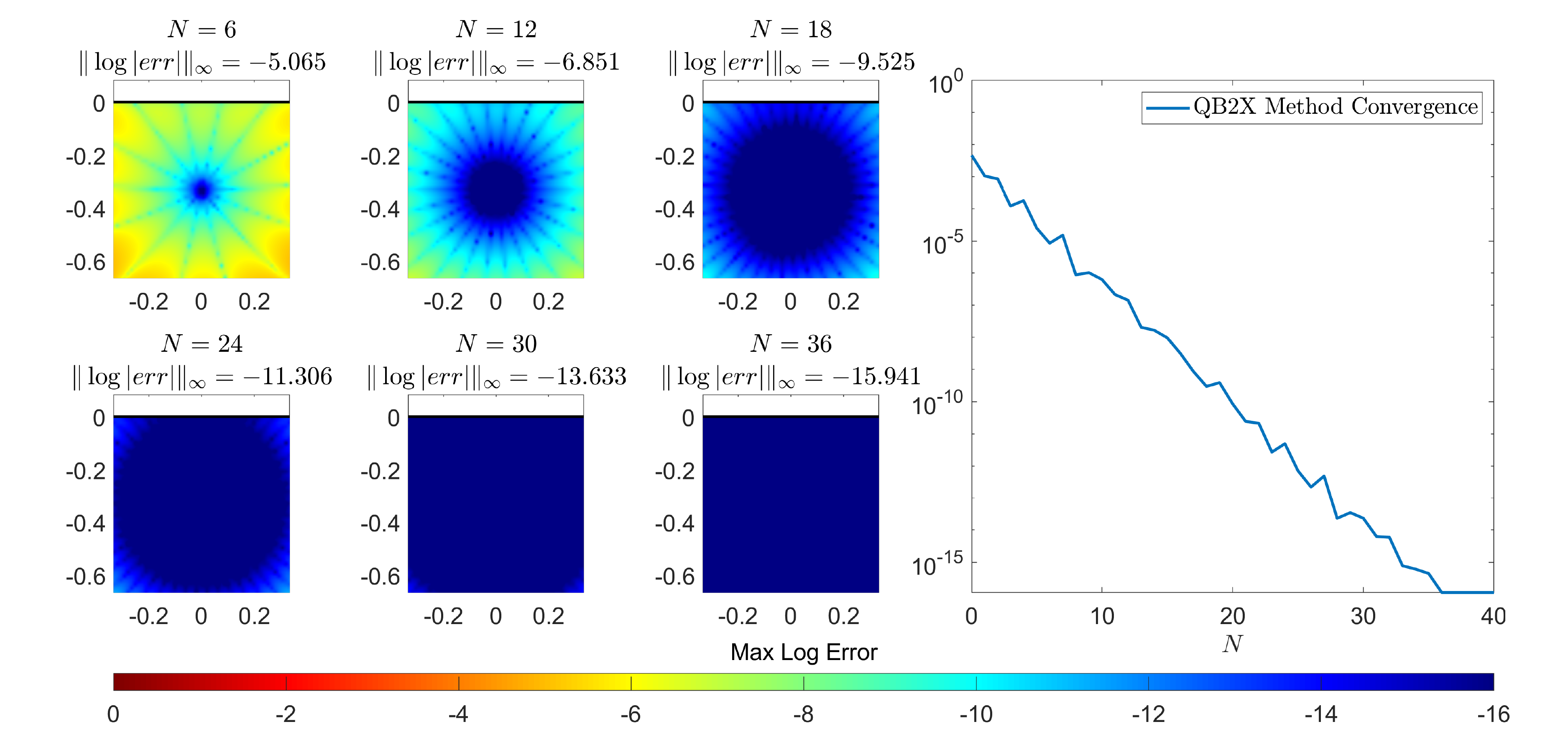} 
   \caption{Double layer potential for $s(t)=0$, $\psi(t)=\frac{2t^2+2t+3}{4}$, and $k=1$: (\textit{Left}) QB2X errors for 
$N=6$, $12$, $18$, $24$, $30$, and $36$. We set $P=50$. (\textit{Right}) Convergence of the $L^\infty$ error for different $N$ values.}
   \label{fig:DLP_s0_rhopoly}
\end{figure}

\begin{figure}[htbp] 
   \centering
   \includegraphics[width=4.5in]{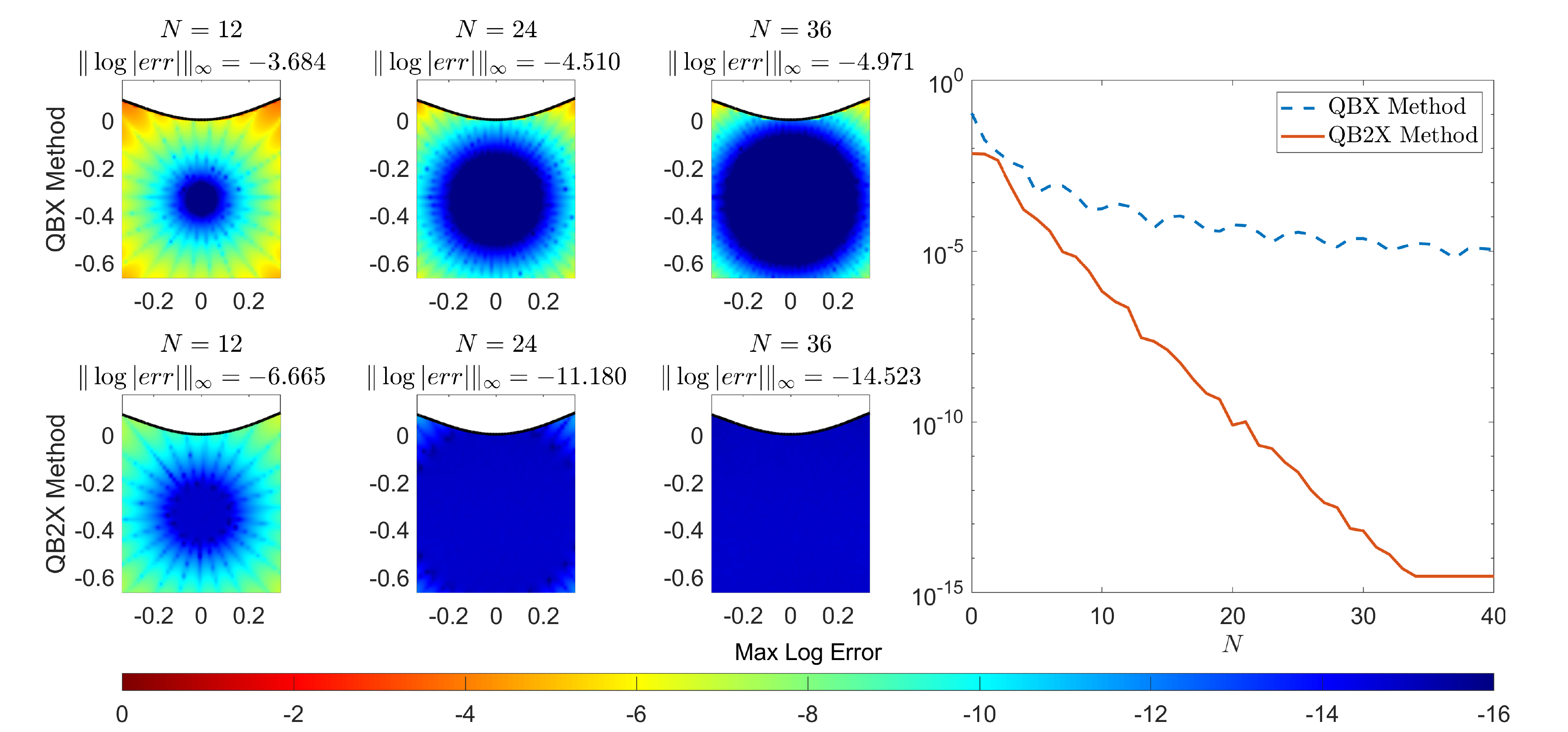} 
	\caption{Double layer potentials for $s(t)=t^2+\frac{t^3}{10}-2t^4$,  $k=1$, and $\psi(t)=1$: (\textit{Left}) QBX (upper) 
	and QB2X (lower) errors for different $N$ values. We set $P=100$. 
  (\textit{Right}) Convergence of the $L^{\infty}$ error for both the QBX (dashed line) and QB2X (solid line) methods for increasing $N$.}
   \label{fig:DLP_smcurve_rhopoly}
\end{figure}

Finally, we show the QB2X performance for the double layer Helmholtz potentials. We found that the numerical 
results for double layer potentials are similar to those for single layer potentials. In Figure \ref{fig:DLP_s0_rhopoly},
we set all the parameters the same as those in Figure \ref{fig:SLP_s0_rhopoly} ($s(t) = 0$, $k=1$, and $\psi(t)=\frac{2t^2+2t+3}{4}$) 
and compute the double layer potential using QB2X-Helmholtz. Both the error distributions and convergence rates are 
similar to the single layer case. In Figure \ref{fig:DLP_smcurve_rhopoly}, we use the same settings as in 
Figure \ref{fig:SLP_smcurve_QBX_QB2X} ($s(t) = t^2+\frac{t^3}{10}-2t^4$, $k=1$, $\psi(t)=1$) and compute
the double layer Helmholtz potential using both the QBX-Helmholtz and QB2X-Helmholtz methods.  Clearly, each
method shows similar error behavior and rate of convergence ($L^{\infty}$ norm) for single and double layer 
potentials.

\section{Conclusions}
\label{sec:4}
In this paper, we present the QB2X method for evaluating the Helmholtz single and double layer potentials in two-dimensional space.
The QB2X-Helmholtz method uses both local complex Taylor expansions and plane wave type expansions to effectively 
capture the nonlinear contributions from the boundary geometry. This method overcomes the convergence issues found in the original 
QBX method. Under minor restrictions on the discretization of the boundary geometry, the QB2X-Helmholtz representation can be valid and
accurate in the entire leaf box of the FMM hierarchical tree structure. Preliminary numerical results confirm the excellent
accuracy and stability properties of QB2X-Helmholtz. The extension of QB2X-Helmholtz to the modified Helmholtz (Yukawa) 
kernel \cite{huang2009fmm} is straight forward and is being implemented. Extending the QB2X method to 3D layer potentials 
is challenging and is currently being studied.


\section*{Statement and Declarations}

\noindent\textbf{Funding} M.H. Cho was supported by NSF grant DMS 2012382 and a grant from the Simons Foundation (No. 404499).  J. Huang was supported by NSF grant DMS 2012451.\\

\noindent\textbf{Competing Interests} The authors have no relevant financial or non-financial interests to disclose.\\

\noindent \textbf{Data Availability} The datasets generated during and/or analysed during the current study are available from the corresponding author on reasonable request.\\

\bibliographystyle{spmpsci}      
\bibliography{QB2Xref}   


\end{document}